\documentclass[11pt]{article}
\usepackage{a4,amsfonts,amssymb,amsmath}
\usepackage{epsfig} 
\usepackage{psfrag}
\usepackage{graphicx}

\def\proof{{\bf Proof:\quad}}
\def\beginpf{\proof}
\def\qed{\hfill\rule{2.2mm}{2.2mm}\vspace{1ex}}
\def\endpf{\qed}

\newtheorem{theorem}{Theorem}[section]

\newtheorem{corollary}[theorem]{Corollary}

\newtheorem{example}[theorem]{Example}

\newtheorem{proposition}[theorem]{Proposition}
\newtheorem{remark}[theorem]{Remark}

\def\eps{\varepsilon}

\def\CC{\mathbb C}

\def\ZZ{\mathbb Z}

\def\NN{\mathbb N}
\def\RR{\mathbb R}

\def\LL{\mathcal L}
\def\JJ{\mathcal J}

\newcommand{\HH}{\mathcal H}

\newcommand{\re}{\mathop{\rm Re}\nolimits}
\renewcommand{\Re}{\re}

\def\text{\mbox}

\title{Applications of Laplace--Carleson embeddings to admissibility and controllability}
\author{{Birgit Jacob\thanks{Fachbereich C - Mathematik und Naturwissenschaften,
Bergische Universit\"at Wuppertal, Gau\ss stra\ss e 20, 42097 Wuppertal, Germany, \tt jacob@math.uni-wuppertal.de}
\qquad Jonathan R.~Partington\thanks{School of Mathematics,
University of Leeds,
Leeds LS2 9JT, U.K. \tt
J.R.Partington@leeds.ac.uk}
\qquad Sandra Pott\thanks{Centre for Mathematics, Faculty of Science, Lund University, S\"olvegatan 18, 22100 Lund, Sweden, \tt Sandra.Pott@math.lu.se}}}

\begin{document}
\maketitle
\begin{abstract}
It is shown how results on Carleson embeddings induced by the Laplace transform can
be use to derive new and more general results concerning the weighted admissibility of control and
observation operators for linear semigroup systems with $q$-Riesz bases of eigenvectors.
Next, a new Carleson embedding result is proved, which gives further results on weighted admissibility
for analytic semigroups.
Finally, controllability by smoother inputs is characterised by means of a new result about
weighted interpolation.
\end{abstract}


{\bf 2000 Subject Classification.} 30D55, 30E05, 47A57, 47D06, 93B05, 93B28.

\section{Introduction}

The main purpose of this note is to show how recent results on Carleson embeddings,
mostly derived in \cite{jpp12},
may be applied to the theory of well-posed linear systems. We shall discuss
concepts such as admissibility, controllability and observability,
for which some standard references are the books \cite{TW,staffans}
and the survey \cite{jp}.
Our basic tools will be the theory of a class of function spaces known as
Zen spaces, which include the standard Hardy and Bergman spaces.\\

The structure of the paper is as follows. In Section \ref{sec:zen}, we review the
basics of the theory of admissibility for diagonal semigroups, introduce the
key embedding results for Zen spaces, and derive new admissibility results in
this context.
Section \ref{sec:pq} treats the case of analytic semigroups (corresponding to
measures supported in a sector), where we are able to move away from the
Hilbertian ($L^2$) context and study general $L^p$ spaces. We also derive
a new embedding theorem and apply it to weighted admissibility.
Finally, in Section \ref{sec:4} we consider notions of controllability, which
are linked with interpolation questions.

\section{Admissibility for diagonal semigroups}  \label{sec:zen}

Let $A$ be the infinitesimal generator of a $C_0$--semigroup $(T(t))_{t \ge 0}$
defined on a Hilbert space $H$, and
consider the system
\begin{equation}\label{eq:ABsystem}
{dx(t) \over dt}=Ax(t)+Bu(t), \qquad x(0)=x_0, \quad t \ge 0,
\end{equation}
where $u(t) \in \mathbb C$ is the input at time $t$,
and 
$B: \mathbb C \to D(A^*)'$, the control operator.
Here $D(A^*)'$ is the completion of $H$ with respect to
the norm
\[
\|x\|_{D(A^*)'}=\|(\beta-A)^{-1}x\|_H,
\]
for any $\beta\in\rho(A)$.
To ensure that the state $x(t)$ lies in $H$
it is sufficient that $B \in \mathcal{L}(\mathbb C,D(A^*)')$ and  
\[
 \left\|\int_0^\infty T(t) Bu(t) \, dt \right\|_H \le m_0 \|u\|_{L^2(0,\infty)}
\]
for some $m_0>0$
(the admissibility condition for $B$). We note that the $C_0$--semigroup $(T(t))_{t \ge 0}$ has an extension to $ D(A^*)'$.

Dually, we may consider the system
\[
{dx(t) \over dt}=Ax(t), \qquad y(t)=Cx(t),
\]
with $x(0)=x_0$, say. 
Here $C:D(A) \to \mathbb C$ is an $A$-bounded observation operator mapping; 
i.e., for some $m_1$, $m_2 > 0$,
$$\|Cz\| \le m_1 \|z\| + m_2\|Az\| .$$

$C$ is {\em admissible}, if there is an $m_0 > 0$ such that
$y \in L^2(0,\infty)$ and
$\|y\|_2 \le m_0 \|x_0\|$ for all $x_0 \in D(A)$.
Note that $y(t)=CT(t) x_0$ for every $x_0\in D(A)$.\\

The duality here is that
$B$ is an admissible control operator for $(T(t))_{t \ge 0}$ if and only if $B^*$ is an admissible 
observation operator for 
the dual semigroup $(T(t)^*)_{t \ge 0}$.\\

Admissibility is an important concept in the theory of well-posed systems, and we refer to the
survey \cite{jp} and the book \cite{TW} for the basic background to the subject.
For diagonal semigroups, admissibility is linked with the theory of Carleson measures as in
\cite{HR83,weiss88}; namely, supposing that $A$ has a Riesz basis of eigenvectors, with
eigenvalues $(\lambda_k)$, then a scalar control operator corresponding to a sequence $(b_k)$
is admissible if and only if the measure
\[
\mu:=\sum_k |b_k|^2 \delta_{-\lambda_k}
\]
is a Carleson measure for the Hardy space $H^2(\CC_+)$ on the right half-plane: this means that
the canonical embedding $H^2(\CC_+) \to L^2(\CC_+,\mu)$ is bounded.
In fact, the extension to normal semigroups
has also been made \cite{weiss99}.

Generalizations to $\alpha$-admissibility, in which
$u$ must lie in $L^2(0,\infty; t^{\alpha} \, dt)$ for $-1<\alpha<0$, were studied by Wynn \cite{wynn}.
The key fact here is that under the Laplace transform the space
$L^2(0,\infty; t^{\alpha} \, dt)$ is mapped to a weighted Bergman space, and for these 
there are analogues of the Carleson measure theorem available.

The results below enable us to take this generalization further and consider admissibility
in the sense of the input lying in much more general spaces $L^2_w(0,\infty)=L^2(0,\infty; w(t) \, dt)$.

In order to state the link between admissibility and embeddings in the greatest generality possible, we assume now
that $1 \le q < \infty$ and the semigroup $(T(t))_{t \ge 0}$ acts on a Banach space $X$
with a $q$-Riesz basis of eigenvectors $(\phi_k)$; that is, $T(t)\phi_k=e^{\lambda_k t}\phi_k$ 
for each $k$, and
$(\phi_k)$ is a
a Schauder basis of $X$ such that
for some $C_1, C_2 > 0$ we have
\[
C_1 \sum |a_k|^q \le \| \sum a_k \phi_k \|^q \le C_2 \sum |a_k|^q
\]
for all sequences $(a_k)$ in $\ell^q$. 
In practice, this will mean that without loss of generality we can assume that $X=\ell^q$ and that
the eigenvectors of the generator of $(T(t))_{t \ge 0}$, denoted by $A$, are the canonical basis of $\ell^q$.
We suppose also that we have a Banach
space $Z$ of functions on $(0,\infty)$, which may be an $L^p$ space or a weighted space
$L^2(0,\infty; w(t)dt)$, whose dual space $Z^*$ can be regarded, respectively, as 
either $L^{p'}(0,\infty)$
or  $L^2(0,\infty; w(t)^{-1}dt)$ in a natural way. Here, and throughout the paper, $p'$
denotes the conjugate exponent to $p$, i.e., $p'=p/(p-1)$. 

The following result links admissibility and Laplace--Carleson embeddings (that is, Carleson embeddings
induced by the Laplace transform).

\begin{theorem}\label{thm:key9}
 Let $B$ be a linear bounded map from $\mathbb C$ to $D(A^*)'$ corresponding to the sequence $(b_k)$. The control operator $B$ is $Z$-admissible for $(T(t))_{t \ge 0}$, that is,
there is a constant $m_0>0$ such that
\[
 \left\|\int_0^\infty T(t) Bu(t) \, dt \right\|_X \le m_0 \|u\|_Z,\quad u\in Z,
\]
if and only if 
the Laplace transform induces a continuous mapping from $Z$
into $L^{q}(\CC_+,d\mu)$, where $\mu$ is the measure $\sum |b_k|^{q}\delta_{-\lambda_k}$.
\end{theorem}
\beginpf 
Clearly we may suppose without loss of generality that $X=\ell^q$ and that $\phi_k=e_k$, the standard
basis of $X$.
We have that 
\begin{eqnarray*}
\left\|\int_0^\infty T(t) Bu(t) \, dt \right\|_X 
&=& \left \| \int_0^\infty \sum_k e^{\lambda_k t} b_k e_k u(t) \, dt \right \|_X \\
&=& \left( \sum_k |\hat u(-\lambda_k)|^q |b_k|^q \right)^{1/q},
\end{eqnarray*}
from which the result follows easily.
\endpf

A duality argument gives the corresponding result for observation operators.

\begin{theorem}\label{thm:key}
 Let $C$ be a linear bounded map from  $D(A)$ to $\mathbb C$. The observation operator $C$ is $Z$-admissible for $(T(t))_{t \ge 0}$, that is,
there is a constant $m_0>0$ such that
\[
\| CT(\cdot)x \|_Z \le m_0 \|x\|_X
\]
for all $x \in D(A)$,
if and only if the Laplace transform induces a continuous mapping from $Z^*$
into $L^{q'}(\CC_+,d\mu)$, where $\mu$ is the measure $\sum |c_k|^{q'}\delta_{-\lambda_k}$ and $c_k:=C \phi_k $.
\end{theorem}

\beginpf
Again we may suppose without loss of generality that $X=\ell^q$ and that $\phi_k=e_k$, the standard
basis of $X$.
$Z$-admissibility is equivalent to the condition that
\[
\sup_{f,x}\left| \int_0^\infty CT(t)x f(t) \, dt \right|< \infty,
\]
where we take $f \in Z^*$ and $x =(x_k)\in D(A)$ both of norm 1.
Calculating $CT(t)x$, we see that this is equivalent to the condition that
\[
\sup_{f,x} | \sum_k c_k \hat f(-\lambda_k) x_k | < \infty,
\]
and, taking the supremum over the set of $x$  of norm 1 in $D(A)$, which is dense, we obtain
\[
\sup_f \|( c_k  \hat f(-\lambda_k)) \|_{q'} < \infty,
\]
which is easily seen to be equivalent to the boundedness of the Laplace--Carleson
embedding from $Z^*$
into $L^{q'}(\CC_+,d\mu)$.
\endpf


In general we shall state our results in terms of control operators, leaving the
interested reader to deduce the corresponding results for observation operators.\\


Now let $\tilde \nu$ be a positive  regular Borel measure on $[0, \infty)$ satisfying the following $(\Delta_2)$-condition:
\begin{equation}\tag{$\Delta_2$}
   R:=  \sup_{r >0} \frac{\tilde \nu[0, 2r)}{\tilde \nu[0, r)} < \infty.
\end{equation}
Let $\nu$ be the positive regular Borel measure on $\overline{\CC_+} =[0, \infty) \times \RR $ given by $d\nu = d \tilde \nu \otimes d \lambda$, where $\lambda$ denotes Lebesgue measure.
In this case, for $1\le p<\infty$, we call 
$$A^p_\nu= \left\{ f: \CC_+ \rightarrow \CC \text{ analytic}: \sup_{\eps >0} \int_{\overline{\CC_+}} |f(z+\eps)|^p d\nu(z)< \infty \right\}
$$ a Zen space on $\CC_+$. If $\tilde \nu(\{0\}) >0$, then by standard Hardy space theory, $f $ has a well-defined boundary function $\tilde f \in L^p(i \RR)$, and we can give meaning
to the expression $\int_{\overline{\CC_+}} |f(z)|^p d\nu(z)$.   Therefore, we write
$$
                         \|f\|_{A^p_\nu} = \left(\int_{\overline{\CC_+}} |f(z)|^p d\nu(z) \right)^{1/p}.
$$
Clearly the space $A^2_\nu$ is a Hilbert space.

Well-known examples of Zen spaces are  Hardy spaces $H^p(\CC_+)$, where $\tilde \nu$ is the Dirac measure at $0$, or the standard weighted Bergman spaces $A^p_{\alpha}$, where 
$d \tilde \nu(r) = r^\alpha dr $, $ \alpha > -1$.  
Some further examples constructed from Hardy spaces on shifted half planes
were given by Zen Harper in \cite{zen09,zen10}.

The following proposition, given in \cite{jpp12}, is elementary and appears for special cases in \cite{zen09,zen10}. Partial results are also given in
\cite{DP94,DGM}.
\begin{proposition} (Proposition 2.3 in \cite{jpp12}) Let $A^2_\nu$ be a Zen space,   and let $w:(0, \infty) \rightarrow \RR_+$ be given by
$$
    w(t) = 2 \pi \int_0^\infty  e^{-2rt} d \tilde \nu(r)   \qquad (t >0).
$$
Then the Laplace transform defines an isometric map $\LL: L^2_w(0, \infty) \rightarrow A^2_\nu$.
\end{proposition}
Note that the existence of the integral is guaranteed by the $(\Delta_2)$-condition.
\bigskip

We shall require the following Laplace--Carleson Embedding Theorem
from \cite{jpp12}.

\begin{theorem} (Theorem 2.4 in \cite{jpp12}) \label{thm:jpp12a} Let $A^2_\nu$ be a Zen space, $ \nu = \tilde \nu \otimes  \lambda$, and let $w:(0, \infty) \rightarrow \RR_+$ be given by
\begin{equation}\label{massw}
    w(t) = 2 \pi \int_0^\infty  e^{-2rt} d \tilde \nu(r)   \qquad (t >0).
\end{equation}
Then the following are equivalent:
\begin{enumerate}
\item The Laplace transform $ \LL$ given by  $ \LL f (z) = \int_0^\infty e^{-t z} f(t) dt$ defines a bounded linear map
$$
\LL: L^2_{w}(0,\infty) \rightarrow L^2(\CC_+, \mu).
$$
\item
For a sufficiently large $N \in \NN$, there exists a constant $\kappa >0$ such that
$$
      \int_{\CC_+}    \left|  (\LL t^{N-1} e^{-\lambda t }) (z)\right|^2 d \mu(z)  
\le \kappa     \int_0^\infty    | t^{N-1} e^{-\lambda t } |^2  w(t) dt         \text{ for each } \lambda \in \CC_+.
$$
 \item There exists a constant $\kappa>0$ such that 
\begin{equation*}   
 \mu(Q_I) \le \kappa \nu(Q_I) \text{ for each Carleson square }Q_I.
 \end{equation*}
where $Q_I$ denotes  the {\em Carleson square} $Q_I = \{ z= x + iy \in \CC_+: iy \in I, 0 < x < |I|\}$.
\end{enumerate}
\end{theorem}

From this we may deduce a result characterizing admissibility for normal semigroups
in the sense of $L^2_{w}(0,\infty)$.

\begin{theorem}\label{Theo:adm1}
Suppose that $A:D(A)\subset H \rightarrow H$ has a Riesz basis $(\phi_k)$ of eigenvectors with eigenvalues $(\lambda_k)$ satisfying Re$\, \lambda_k<0$ and let $B$ be a linear bounded map from $\mathbb C$ to $D(A^*)'$ given by the sequence $(b_k)$. Let $A^2_\nu$ be a Zen space and let $w$ be given by (\ref{massw}). Then the following statements are equivalent.
\begin{enumerate}
\item$B$ is an admissible control operator with respect to $L^2_{w}(0,\infty)$, that is, there exists a constant $m_0>0$ such that
\[
 \left\|\int_0^\infty T(t) Bu(t) \, dt \right\|_H \le m_0 \|u\|_{L^2_w(0,\infty)}
\]
for every $u\in L^2_w(0,\infty)$.
\item For a sufficiently large $N \in \NN$, there exists a constant $\kappa >0$ such that
\[
\|(\lambda-A)^{-N}B\|^2 \le \kappa\int_0^\infty |t^{N-1}e^{-\lambda t}|^2 w(t) \, dt, \qquad (\lambda \in \CC_+).
\]
\item There exists a constant $\kappa>0$ such that 
\begin{equation*}   
 \mu(Q_I) \le \kappa \nu(Q_I) \text{ for each Carleson square }Q_I,
 \end{equation*}
where $\mu=\sum_k |b_k|^2 \delta_{-\lambda_k}$.
\end{enumerate}
\end{theorem}

\beginpf
This follows from Theorem \ref{thm:key9} and Theorem \ref{thm:jpp12a}, taking $q=2$
and $Z=L^2_{w}(0,\infty)$. Note that the resolvent condition follows because
\[
\|(\lambda - A)^{-N}B\|^2 = \sum_k |\lambda-\lambda_k|^{-2N}|b_k|^2 = \int_{\CC_+} \frac{d\mu(z)}{|\lambda+z|^{2N}},
\]
and the Laplace transform of $t^{N-1}e^{-\lambda t}$ is a constant multiple of $(\lambda+z)^{-N}$.
\endpf

\begin{remark}\label{rem:classics}
Theorem \ref{Theo:adm1} in the case that $\tilde\nu$ equals the Dirac measure in $0$, that is $A_\nu^2=H^2(\mathbb C_+)$, is due to Ho and Russell \cite{HR83}, and Weiss \cite{weiss88}.
In the case $d\tilde\nu=r^\alpha dr$, $\alpha\in(-1,0)$ the result is due to Wynn \cite{wynn}. Using Theorem \ref{Theo:adm1} any
$\alpha > -1$ can now be considered. 
\end{remark}

In \cite{haak}, Haak applied Carleson
measure theory to
find conditions for admissibility of the system (\ref{eq:ABsystem}) above, with $A$ generating a diagonal semigroup defined on $\ell^q$, and inputs 
lying in the space
$L^p(0,\infty)$. Using the Hausdorff--Young inequality, it is possible to  characterize admissibility in the situation  $p\le 2$ and $1<p'\le q<\infty$. 

\begin{theorem}\label{Theo:adm2}
Let $p\le 2$ and  $1<p'\le q<\infty$. Suppose that $A:D(A)\subset \ell^q \rightarrow \ell^q$ is a diagonal operator with eigenvalues  $(\lambda_k)$ satisfying Re$\, \lambda_k<0$, and let $B$ be a linear bounded map from $\mathbb C$ to $D(A^*)'$ corresponding to the sequence $(b_k)$.  Then the following statements are equivalent.
\begin{enumerate}
\item $B$ is an admissible control operator with respect to $L^{p}(0,\infty)$, that is, there exists a constant $m_0>0$ such that
\[
 \left\|\int_0^\infty T(t) Bu(t) \, dt \right\|_{\ell^q} \le m_0 \|u\|_{L^p(0,\infty)}
\]
for every $u\in L^p(0,\infty)$.
\item There exists a constant $\kappa>0$ such that 
\begin{equation*}   
 \mu(Q_I) \le \kappa |I|^{q/p'} 
 \end{equation*}
for all intervals in $I \subset i\RR$,
where $\mu=\sum_k |b_k|^{q} \delta_{-\lambda_k}$.
\item There exists a constant $\kappa>0$ such that 
\[ \|(\lambda-A)^{-1}B\|_{\ell^q} \le \kappa |\mbox{\rm Re}\, \lambda|^{-1/p}\]
for all $\lambda \in \mathbb C_+$.
\end{enumerate}
\end{theorem}

The theorem is a corollary of Theorem \ref{thm:key9}
combined with \cite[Thm.~3.2]{jpp12}. Moreover the equivalence of Part 1 and 2 can be found in \cite{haak}.

\section{Admissibility for analytic semigroups}    \label{sec:pq}

There is no general characterization of admissibility of the system (\ref{eq:ABsystem}) above, with $A$ generating a diagonal semigroup defined on $\ell^q$, and inputs 
lying in the space $L^p(0,\infty)$, $p>2$, as 
there is no known full characterization of boundedness
of Laplace--Carleson embeddings 
$$
  L^p(0,\infty) \rightarrow L^{q}(\CC_+, \mu), \quad f
  \mapsto \LL f= \int_0^\infty e^{-t \cdot} f(t) dt.
$$
However, characterizations are possible in some cases with additional information on the support on the measure.

 If the measure $\mu$ is supported on a sector $S(\theta) = \{
z \in \CC_+: |\arg z| < \theta \}$ for some $0< \theta <
\frac{\pi}{2}$,
then the oscillatory part of the Laplace transform can be discounted,
and a full characterization of boundedness can be achieved (see also
\cite{haak}, Theorem 3.2 for an alternative characterization by means
of a different measure).

\begin{theorem} (Theorem 3.3 in \cite{jpp12})\label{thm:pqemb} Let $\mu$ be a positive regular Borel measure supported in a
  sector $S(\theta) \subset \CC_+$, $0 < \theta < \frac{\pi}{2}$, and
  let $q \ge p >1$. Then the following are equivalent:
\begin{enumerate}
\item The Laplace--Carleson embedding
$$
     \LL: L^p(0,\infty) \rightarrow L^q(\CC_+, \mu), \quad f \mapsto \LL f,
$$
is well-defined and bounded.

\item There exists a constant $\kappa>0$ such that $\mu(Q_I) \le \kappa
  |I|^{q/p'}$ for all intervals in $I \subset i\RR$ which are symmetric
  about $0$.

\item  There exists a constant $\kappa>0$ such that $\|\LL e^{- \cdot z}
  \|_{L^q_\mu} \le \kappa \|e^{- \cdot z}\|_{L^p}$ for all $z \in \RR_+$.
  
\end{enumerate}
\end{theorem} 

From this we may deduce a result characterizing admissibility for analytic semigroups with respect to $L^p(0,\infty)$.

\begin{theorem}\label{Theo:adm2a}
Let $1<p\le q<\infty$. Suppose that $A:D(A)\subset \ell^q \rightarrow \ell^q$ is a diagonal operator with eigenvalues  $(\lambda_k)$ satisfying Re$\, \lambda_k<0$ and $(-\lambda_k)\subset S(\theta)$ for some $\theta \in (0,\frac{\pi}{2})$, and let $B$ be a linear bounded map from $\mathbb C$ to $D(A^*)'$
given by the sequence $(b_k)$.  Then the following statements are equivalent.
\begin{enumerate}
\item $B$ is an admissible control operator with respect to $L^{p}(0,\infty)$, that is, there exists a constant $m_0>0$ such that
\[
 \left\|\int_0^\infty T(t) Bu(t) \, dt \right\|_{\ell^q} \le m_0 \|u\|_{L^p(0,\infty)}
\]
for every $u\in L^p(0,\infty)$.
\item There exists a constant $\kappa>0$ such that 
\begin{equation*}   
 \mu(Q_I) \le \kappa \nu|I|^{q/p'} 
 \end{equation*}
for all intervals in $I \subset i\RR$ which are symmetric
  about $0$,
where $\mu=\sum_k |b_k|^{q} \delta_{-\lambda_k}$.
\item There exists a constant $\kappa>0$ such that 
\[ \|(z-A)^{-1}B\|_{\ell^q} \le \kappa z^{-1/p}\]
for all $z \in \mathbb R_+$.
\end{enumerate}
\end{theorem}

\begin{remark}
Let $\mu$, $\theta$, $p$ and $q$ be as in Theorem \ref{Theo:adm2a}. 
In \cite{haak}, Theorem 3.2, the equivalence of the following statements is shown:
\begin{enumerate}
\item $B$ is an admissible control operator with respect to $L^{p}(0,\infty)$, that is, there exists a constant $m_0>0$ such that
\[
 \left\|\int_0^\infty T(t) Bu(t) \, dt \right\|_{\ell^q} \le m_0 \|u\|_{L^p(0,\infty)}
\]
for every $u\in L^p(0,\infty)$.
\item There exists a constant $\kappa>0$ such that $\tilde \mu(Q_I) \le \kappa
  |I|^{q/p}$ for all intervals in $I \subset i\RR$ which are symmetric
  about $0$, where $d\tilde \mu(z) = |z|^{q} d\mu(\frac{1}{z})$.
  \end{enumerate}
\end{remark}


Now let us consider the case $p>q$ for sectorial measures $\mu$. In \cite{jpp12} a condition in terms of the balayage $S_{ \mu}$ of $ \mu$ has been obtained.
 Recall that the balayage $S_\mu$ of a positive Borel measure $\mu$ on $\CC_+$ is given by
$S_\mu(t)=\int_{\CC_+} p_z(t) d \mu(z)$, where $p_z$ denotes the Poisson kernel of the right half plane. Let $S_n:= \{z\in\mathbb C \mid 2^{n-1}<{\rm Re}\, z \le 2^n\}$.

\begin{theorem} (Theorem 3.5 in \cite{jpp12}) \label{thm:pgqemb}
 Let $\mu$ be a positive regular Borel measure supported in a
  sector $S(\theta) \subset \CC_+$, $0 < \theta < \frac{\pi}{2}$ and let $1 \le q  <
 p < \infty$. Then the following are equivalent:
 \begin{enumerate}
 \item    The embedding
$$
     \LL: L^p(\RR_+) \rightarrow L^q(\CC_+, \mu), \quad f \mapsto \LL f,
$$
is well-defined and bounded.
\item The sequence $(2^{-n q/p'} \mu(S_{n}))$ is in $\ell^{p/(p-q)}(\ZZ)$.

\item The sequence $(       2^{n/p} \| \LL e^{-2^n}\|_{L^q_\mu} ) $   is in $\ell^{qp/(p-q)}(\ZZ)$.
\end{enumerate}
If $p'<q$, then the above is also equivalent to
\begin{enumerate}
\setcounter{enumi}{3}
\item   $t^{q(2-p)/p}S_{\mu} \in L^{p/(p-q)}(\RR)$.
\end{enumerate}
\end{theorem}

As a corollary we obtain the following result.

\begin{theorem}\label{Theo:adm2b}
Let $1\le q < p<\infty$. Suppose that $A:D(A)\subset \ell^q \rightarrow \ell^q$ is a diagonal operator with eigenvalues  $(\lambda_k)$ satisfying Re$\, \lambda_k<0$ and $(-\lambda_k)\subset S(\theta)$ for some $\theta \in (0,\frac{\pi}{2})$, and let $B$ be a linear bounded map from $\mathbb C$ to $D(A^*)'$
given by the sequence $(b_k)$.  Then the following statements are equivalent.
\begin{enumerate}
\item $B$ is an admissible control operator with respect to $L^{p}(0,\infty)$, that is, there exists a constant $m_0>0$ such that
\[
 \left\|\int_0^\infty T(t) Bu(t) \, dt \right\|_{\ell^q} \le m_0 \|u\|_{L^p(0,\infty)}
\]
for every $u\in L^p(0,\infty)$.
\item The sequence $(2^{-n q/p'} \mu(S_{n}))$ is in $\ell^{p/(p-q)}(\ZZ)$.

\item The sequence $(  2^{n/p} \| (2^n-A)^{-1}B\|_{\ell^{q}} ) $   is in $\ell^{qp/(p-q)}(\ZZ)$.
\end{enumerate}
If $p'<q$, then the above is also equivalent to
\begin{enumerate}
\setcounter{enumi}{3}
\item   $t^{q(2-p)/p}S_{\mu} \in L^{p/(p-q)}(\RR)$.
\end{enumerate}
Here $\mu=\sum_k |b_k|^{q} \delta_{-\lambda_k}$.
\end{theorem}

\begin{example}
We study the one-dimensional heat equation on the interval $[0,1]$ which is given by
\begin{eqnarray*}
\frac{\partial z}{\partial t} (\zeta ,t)&=& \frac{\partial^2 z}{\partial \zeta^2} (\zeta ,t),\qquad \zeta\in(0,1), t\ge 0,\\
\frac{\partial z}{\partial \zeta} (0 ,t) &=& 0, \quad  \frac{\partial z}{\partial \zeta} (1,t)=u(t),\quad t\ge 0,\\
z(\zeta,0) &=& z_0(\zeta), \quad \zeta\in(0,1).
\end{eqnarray*}
This PDE can be written equivalently in the form (\ref{eq:ABsystem}) with $X=\ell^2$, $A e_n = -n^2\pi^2  e_n $, and $b_n$ defined by $b_n=1$
for each $n$. By Theorem \ref{Theo:adm2a} (for $1 \le p \le 2$) and Theorem \ref{Theo:adm2b}
(for $2 \le p < \infty$), the operator $B$ is an admissible control operator with respect to $L^p(0,\infty)$ if and only if $p \ge 4/3$.
\end{example}


Control problems involving smoother, Sobolev--space valued, controls are related to  embeddings of the
form
$$
  \HH^p_{\beta, w}(0,\infty) \rightarrow L^q(\CC_+, \mu), \quad f
  \mapsto \LL f= \int_0^\infty e^{-t \cdot} f(t) dt,
$$
given by the Laplace transform $\LL$. Here,  
for $\beta>0$ the space  $\HH^p_\beta(0,\infty)$ is given by
\begin{eqnarray*}
 \HH^p_{\beta}(0,\infty) &=& \left\{ f \in L^p(\RR_+): \int_0^\infty
    \left| \left( \frac{d}{dx} \right)^\beta f(t) \right|^p dt < \infty\right\},\\
    \|f\|_{\HH^p_\beta}^p &=&\|f\|^p_p + \left\| \left(\frac{d}{dx} \right)^\beta f \right\|^p_p.
\end{eqnarray*}
Here $(\frac{d}{dx})^\beta f$ is defined as a fractional derivative
via the Fourier transform. 
We call $\HH^p_{\beta, w}$
the Sobolev space of index $\beta$ and weight $w$.

In \cite[Corollaries 3.7 and 3.8]{jpp12} the following characterisations of the boundedness of the Laplace--Carleson embedding has been proved.

\begin{proposition} (Corollary 3.7 in \cite{jpp12}) \label{thm:pqsobemb} Let $\mu$ be a positive Borel measure supported in a
  sector $S(\theta) \subset \CC_+$, $0 < \theta < \frac{\pi}{2}$, and
  let $q \ge p >1$. Then the following are equivalent:
\begin{enumerate}
\item The Laplace--Carleson embedding
$$
     \LL: \HH_\beta^p(0,\infty) \rightarrow L^q(\CC_+, \mu), \quad f \mapsto \LL f,
$$
is well-defined and bounded.

\item There exists a constant $\kappa>0$ such that $\mu_{q,\beta}(Q_I) \le \kappa
  |I|^{q/p'}$ for all intervals in $I \subset i\RR$ which are symmetric
  about $0$. Here, $d\mu_{q,\beta}(z) = (1+ \frac{1}{|z|^{q \beta}}) d \mu(z)$.

\item  There exists a constant $\kappa>0$ such that $\|\LL e^{- \cdot z}
  \|_{L^q_\mu} \le \kappa \|e^{- \cdot z}\|_{\HH_\beta^p}$ for all $z \in \RR_+$.

\end{enumerate}
\end{proposition}

\begin{proposition} (Corollary 3.8 in \cite{jpp12}) \label{thm:pgqsobemb}
 Let $\mu$ be a positive regular Borel measure supported in a
  sector $S(\theta) \subset \CC_+$, $0 < \theta < \frac{\pi}{2}$ and let $1 \le q  <
 p$, $\beta \ge 0$. Suppose that $S_{ \mu_{\beta,q}} \in L^{p/(p-q)}$. 
 Then the embedding
$$
     \LL: \HH_{\beta}^p(0,\infty) \rightarrow L^q(\CC_+, \mu), \quad f \mapsto \LL f,
$$
is well-defined and bounded.
\end{proposition}

As an application we are able to characterize admissibility with respect to Sobolev space valued control functions.

\begin{theorem}\label{Theo:adm3}
Let $1<p < q<\infty$. Suppose that $A:D(A)\subset \ell^q \rightarrow \ell^q$ is a diagonal operator with eigenvalues  $(\lambda_k)$ satisfying Re$\, \lambda_k<0$ and $(-\lambda_k)\subset S(\theta)$ for some $\theta \in (0,\frac{\pi}{2})$, and let $B$ be a linear bounded map from $\mathbb C$ to $D(A^*)'$.  Write  $\mu=\sum |b_k|^q \delta_{-\lambda_k}$. Then the following statements are equivalent.
\begin{enumerate}
\item $B$ is an admissible control operator with respect to $\HH_\beta^p(0,\infty)$, that is, there exists a constant $m_0>0$ such that
\[
 \left\|\int_0^\infty T(t) Bu(t) \, dt \right\|_{\ell^q} \le m_0 \|u\|_{\HH_\beta^p(0,\infty)},
\]
for every $u\in \HH_\beta^p(0,\infty)$.
\item There exists a constant $\kappa>0$ such that $\mu_{q,\beta}(Q_I) \le \kappa
  |I|^{q/p'}$ for all intervals in $I \subset i\RR$ which are symmetric
  about $0$. Here, $d\mu_{q,\beta}(z) = (1+ \frac{1}{|z|^{q \beta}}) d \mu(z)$.
\item  There exists a constant $\kappa>0$ such that $\|(z-A)^{-1}B
  \|_{\ell^q} \le \kappa \|e^{- \cdot z}\|_{\HH_\beta^p}$ for all $z \in \RR_+$.

\end{enumerate}
\end{theorem}

\begin{theorem}\label{Theo:adm4}
Let $1<q< p<\infty$. Suppose that $A:D(A)\subset \ell^q \rightarrow \ell^q$ is a diagonal operator with eigenvalues  $(\lambda_k)$ satisfying Re$\, \lambda_k<0$ and $(-\lambda_k)\subset S(\theta)$ for some $\theta \in (0,\frac{\pi}{2})$, and let $B$ be a linear bounded map from $\mathbb C$ to $D(A^*)'$.  Suppose that $S_{\tilde \mu_{\beta,q}} \in L^{p/(p-q)}$. Then  $B$ is an admissible control operator with respect to $\HH_\beta^p(0,\infty)$, that is, there exists a constant $m_0>0$ such that
\[
 \left\|\int_0^\infty T(t) Bu(t) \, dt \right\|_{\ell^q} \le m_0 \|u\|_{\HH_\beta^p(0,\infty)},
\]
for every $u\in \HH_\beta^p(0,\infty)$.
\end{theorem}

In the setting of analytic diagonal semigroups, it is also possible to characterise 
$L^2((0,\infty),t^{\alpha}dt)$-admissibility of control operators and $L^2((0,\infty),t^{-\alpha}dt)$-admissibility of observation operators
 in the range $0 < \alpha<1$ in terms of a Carleson-type condition or a resolvent condition 
(compare this with the counterexample by Wynn in \cite{wynn1} for a diagonal semigroup).

Generally, the difficulty in the case $\alpha >0$ stems from the well-known fact that 
the boundedness of Carleson embeddings on Dirichlet space cannot be tested on reproducing kernels or by means of a simple Carleson-type condition
\cite{steg}. In the 
case of sectorial measures, however, such a characterisation is possible, at least for
$0 < \alpha < 1$. For a interval $I \subset i \RR$, let $T_I$ denote the right half of the Carleson square $Q_I$.

\begin{theorem} \label{thm:diriemb} Let $\mu$ be a positive Borel measure supported in a
  sector $S(\theta) \subset \CC_+$, $0 < \theta < \frac{\pi}{2}$, and
  let $0<\alpha <1$. Then the following are equivalent:
\begin{enumerate}
\item The Laplace--Carleson embedding
$$
     \LL: L^2((0, \infty), t^{\alpha} dt) \rightarrow L^2(\CC_+, \mu), \quad f \mapsto \LL f,
$$
is well-defined and bounded.

\item There exists a constant $\gamma>0$ such that 
$$\mu(T_I) \le \gamma
  |I|^{1- \alpha}
  $$ 
 for all intervals in $I \subset i\RR$ which are symmetric
  about $0$. 

\item  There exists a constant $\kappa>0$ such that 
\[
\|\LL t^{-\alpha} e^{- t z}
  \|_{L^2(\CC_+,\mu)} \le \kappa \|t^{-\alpha} e^{- t z}\|_{L^2(t^{\alpha}dt)}
\]
 for all $z \in \RR_+$.

\end{enumerate}
\end{theorem}
\proof The implication (1) $\Rightarrow$ (3) is immediate. For  (3) $\Rightarrow$ (2), let  $z=z_I=|I|/2$ denote the centre of the Carleson square $Q_I$ over
an intervals $I \subset i\RR$ which is symmetric
  about $0$. Then the modulus of the function
\begin{equation}\label{eq:nastylap}
     (\LL t^{-\alpha} e^{- t z_I}) (s) = \frac{\Gamma(1-\alpha) }{(z_I +s)^{-\alpha+1}}
\end{equation}
is bounded below by $ \Gamma(1-\alpha) \frac{1}{(2 z_I)^{-\alpha+1}}$ on $T_I$, and therefore
\begin{eqnarray*}
      \mu(T_I) &\le& \frac{(2 z_I)^{2-2\alpha}    }{\Gamma(1-\alpha)^2}   \int_{\CC_+}     |(\LL t^{-\alpha} e^{- t z_I}) (s)|^2   d\mu(s)   \\   
&\le &  \kappa^2 \frac{(2 z_I)^{2-2\alpha}    }{\Gamma(1-\alpha)^2}  \|t^{-\alpha} e^{- t z_I}\|^2_{L^2(t^{\alpha}dt)} 
         = \frac{\kappa^2  }{\Gamma(1-\alpha)} |I|^{1-\alpha}.
\end{eqnarray*}

Let us now consider (2) $\Rightarrow$ (1). We use the argument from \cite{jpp12}, Thm 3.3.
For $n \in \ZZ$, let 
\[
T_n = \{ x+ iy \in \CC_+: 2^{n-1} < x \le 2^{n}, -2^{n-1}< y  \le 2^{n-1} \}.
\]
  That is, $T_n$ is the right half of the Carleson square $Q_{I_n}$
  over the interval $I_n = \{ y \in i\RR, |y| \le 2^{n-1} \}$. The
  $T_n$ are obviously pairwise disjoint.

Without loss of generality we assume  $0<\theta < \arctan(\frac{1}{2}) $, in which
case $S(\theta) \subseteq \bigcup_{n = - \infty}^\infty T_n$ (otherwise, we also have to use finitely many translates $T_{n,k}$ of each $T_n$, for which the same estimates apply) .

Now let $z \in T_n$ for some $n \in \ZZ$. Then we obtain, for $f\in L^2((0,\infty), t^{\alpha} dt)$,
\begin{multline*}
 |\LL f(z)| \le \int_0^\infty |e^{-zt}| |f(t) | dt \le \int_0^\infty
 |t^{-\alpha/2}  e^{-2^{n-1}t}| |f(t)    t^{\alpha/2}       | dt \\  \le C_\alpha 2^{(-n+1)(1-\alpha/2)} (M(t^{\alpha/2}f))(2^{-n+1}),
\end{multline*}
where $C_\alpha >0$ is a constant dependent only on the $L^1$-integration kernel $\phi_{\alpha}(t)= \chi_{[0,\infty)}(t+1)(t+1)^{\alpha/2}  e^{-t-1}$, and 
$Mf$ is the Hardy--Littlewood maximal function. We refer to e.g. \cite{stein}, page 57, equation (16) for a pointwise estimate between the maximal
function induced by the kernel $\phi_\alpha$, and $M$. Note that we can easily dominate $\phi_\alpha$ by a positive, radial, decreasing $L^1$ function here. 
Consequently,
\begin{eqnarray*}
     \int_{S(\theta)} |\LL f(z)|^2 d\mu(z)  
&\le& C_\alpha^2 \sum_{n= -
     \infty}^\infty 2^{(-n+1)(2 - \alpha)} (M (t^{\alpha/2}f)(2^{-n+1}))^2 \mu(T_n)\\
& \le &  \gamma C_\alpha^2 \sum_{n= -
     \infty}^\infty 2^{(-n+1)(2-\alpha)} 2^{n(1-\alpha)} (M(t^{\alpha/2}f)(2^{-n+1}))^2  \\
    &  = & \gamma C_\alpha^2 2^{1-\alpha} \sum_{n= -
     \infty}^\infty  2^{-n+1} M(t^{\alpha/2}f)(2^{-n+1})^2 \\
& \lesssim & \|f\|^2_{L^2(t^{\alpha}dt)}.
\end{eqnarray*}
\qed

We therefore have the following corollary of Theorem  \ref{thm:diriemb} and Theorem \ref{thm:key9}.

\begin{corollary} 
Let $0<\alpha<1$. Suppose that $A:D(A)\subset \ell^2 \rightarrow \ell^2$ is a diagonal operator with eigenvalues  $(\lambda_k)$ satisfying Re$\, \lambda_k<0$ and $(-\lambda_k)\subset S(\theta)$ for some $\theta \in (0,\frac{\pi}{2})$, and let $B$ be a linear bounded map from $\mathbb C$ to $D(A^*)'$.  Write  $\mu=\sum |b_k|^2 \delta_{-\lambda_k}$. Then the following statements are equivalent.
\begin{enumerate}
\item $B$ is an admissible control operator with respect to $L^2((0,\infty),t^{\alpha}dt)$, that is, there exists a constant $m_0>0$ such that
\[
 \left\|\int_0^\infty T(t) Bu(t) \, dt \right\|_{\ell^2} \le m_0 \|u\|_{L^2((0,\infty),t^{\alpha}dt)},
\]
for every $u\in L^2((0,\infty),t^{\alpha}dt)$.
\item There exists a constant $\kappa>0$ such that $\mu (T_I) \le \kappa
  |I|^{1-\alpha}$ for all intervals in $I \subset i\RR$ which are symmetric
  about $0$. 
\item \label{cond:9993} There exists a constant $\kappa>0$ such that 
\[
\|\LL t^{-\alpha}e^{- \lambda t}
  \|_{L^2(\CC_+,\mu)} \le \kappa \|t^{-\alpha}e^{- \lambda t}\|_{L^2(t^\alpha dt)}
\]
 for all $\lambda \in \RR_+$.
\item \label{cond:9994} There exists a constant $\kappa>0$ such that
\[
\|(\lambda-A)^{\alpha-1} B \| \le \kappa \lambda^{(\alpha-1)/2} 
\]
for all $\lambda \in \RR_+$.
 \end{enumerate}
\end{corollary}

\beginpf
It is enough to check that conditions \ref{cond:9993} and \ref{cond:9994} are equivalent.
Note that  for $\lambda \in \RR_+$
\[
\|(\lambda-A)^{\alpha-1}B\|^2 = \sum_k |b_k|^2|\lambda-\lambda_k|^{2\alpha-2},
\]
which is a constant multiple of $\|\LL t^{-\alpha}e^{- \lambda t}
  \|_{L^2(\CC_+,\mu)}^2$ (cf. (\ref{eq:nastylap})) . Similarly, we have that $ \|t^{-\alpha}e^{- \lambda t}\|_{L^2(t^\alpha dt)}^2$ is a constant multiple
of $\lambda^{\alpha-1}$.
\endpf

Note that the corollary is also valid  in the limiting case $\alpha=0$.


\section{Exact controllability for diagonal systems}
\label{sec:4}


We consider the equation
\begin{equation}\label{eqn:cont}
{dx(t) \over dt}=Ax(t)+Bu(t), 
\end{equation}
with  solution
$$ x(t)=T(t) x_0 + \int_0^t T(t-s) Bu(s) \, ds,$$
suitably interpreted.
Again we consider an exponentially stable $C_0$ semigroup $(T(t))_{t \ge 0}$, on a Hilbert space $H$,
i.e., 
$$
\|T(t)\| \le Me^{-\lambda t}, \qquad (t \ge 0),
$$
for some $M>0$ and $\lambda>0$.
Suppose first that $B$ is admissible.
Then we have a bounded operator 
$\mathcal{B}_\infty: L^2(0,\infty;U) \to H$, defined by
\[
\mathcal{B}_\infty u = \int_0^\infty T(t) B u(t) \, dt.
\]
The system is {\em exactly controllable}, if its range
$R( \mathcal{B}_\infty)$ equals $ H$. If $B$ is not admissible, then the operator $\mathcal{B}_\infty$
is commonly defined as a mapping into a larger (extrapolation) space and exact controllability 
requires that its image contains $H$.

In \cite{jp06} exact controllability for diagonal systems with scalar inputs was characterised in terms of Carleson
measures (a version for multivariable inputs was given in \cite{jpp07}). In particular if $A$ has a Riesz basis of eigenvectors,
with eigenvalues $(\lambda_n)$, then with a control operator corresponding to a sequence $(b_n)$ 
the system is exactly controllable if and only if
\[
\nu_\lambda:=\sum_n  \frac{|\re \lambda_n|^2}{|b_n|^2 \prod_{k \ne n} p(\lambda_n,\lambda_k)^2}\delta_{-\lambda_n}
\]
is a Carleson measure. Here $p(\lambda_n,\lambda_k)$ is the pseudo-hyperbolic metric, i.e., 
\[
p(\lambda_n,\lambda_k) = \left |\frac{\lambda_n-\lambda_k}{\lambda_n+\overline \lambda_k}
\right| .
\]
Exact controllability by inputs in Sobolev spaces $\mathcal{H}^2_\beta$ with $0< \beta<1/2$ was characterised in \cite[Thm.~3.8]{JPP09}.
In \cite[Thm.~3.9]{jpp12} the following result was proved, which enable us to dispense with the restriction on $\beta$.


\begin{theorem}
 \label{thm:sobolevcarl}
Let $\mu$ be a positive Borel measure on the  right half plane
$\CC_+$ and let $\beta >0$. Then the following are equivalent:
\begin{enumerate}
\item The Laplace--Carleson embedding
$$
     \HH^2_\beta(0,\infty) \rightarrow L^2( \CC_+, \mu)
$$
is bounded.
\item The measure 
$|1+ z|^{-2 \beta} d \mu(z)$ is a Carleson measure on $\CC_+$.
\end{enumerate}
\end{theorem}

In \cite{ss}, \cite{mcphail}, the following theorem was  proved:

\begin{theorem} \label{thm:weightedint} Let
$(g_k)_{k \in \NN}$ be a sequence of nonzero complex numbers and let
$(z_k)_{k \in \NN}$ be a Blaschke sequence in $\CC_+$.
Write $b_{\infty,k}=\prod_{j\ne k} p(z_j,z_k)$.

Let $m_{2}= \sup_{(a_k) \in \ell^2, \|(a_k)\|_2=1}  \inf_{f \in H^2(\CC_+), g_k f(z_k) = a_k} \|f\|_{H^2}$.

Then
$$
 m_{2} = \|  \JJ_{\mu_{2}}\|
$$
where
$$
   \mu_{2} = \sum_{k=1}^\infty
     \frac{|2\Re  z_k|^{2}}{|b_{\infty,k} g_k|^{2}} \delta_{z_k},
$$
and $\JJ_{\mu_{2}}$ is the Carleson embedding
$$
  \JJ_{\mu_{2}}: H^{2}(\CC_+) \rightarrow L^{2}(\CC_+,\mu_{2}).
$$
\end{theorem}

Using  Theorem \ref{thm:sobolevcarl}, we obtain
\begin{corollary}\label{cor1}
Let $(g_k)_{k \in \NN}$ be a sequence of nonzero complex numbers and let
$(z_k)_{k \in \NN}$ be a Blaschke sequence in $\CC_+$. 
Write $b_{\infty,k}=\prod_{j\ne k} p(z_j,z_k)$.
Let 
$$
   m_\beta = \sup_{(a_k) \in \ell^2, \|(a_k) =1\|} \inf_{f \in
   \HH^2_\beta} \{ \|f\|_{\HH^2_\beta}: g_k \LL f(z_k) =a_k \}.
$$
Then $m_\beta < \infty$, if and only if there is a constant $\kappa>0$ such
that
$$
    \sum_{z_k \in Q_I}  \frac{|2\re  z_k|^{2}|1+ z_k |^{2\beta}}{|b_{\infty,k}|^2
    |g_k|^{2} } \le \kappa |I|     \text{ for all
    intervals } I \subset i\RR.
$$
\end{corollary}

Finally, Corollary \ref{cor1} enables us to characterize controllability by inputs in Sobolev spaces $\mathcal{H}^2_\beta$.

\begin{theorem}\label{Theo:contr1}
Let $\beta>0$. Suppose that $A:D(A)\subset H \rightarrow H$ has a Riesz basis $(\phi_k)$ of eigenvectors with eigenvalues $(\lambda_k)$ satisfying Re$\, \lambda_k<0$ and let $B$ be a linear bounded map from $\mathbb C$ to $D(A^*)'$. 
Write $b_{\infty,k}=\prod_{j\ne k} p(\lambda_j,\lambda_k)$.
Then the following statements are equivalent.
\begin{enumerate}
\item System (\ref{eqn:cont}) is exactly controllable with respect to  $\mathcal{H}^2_\beta$, that is,
\[ H\subset {\cal B}_\infty ( \mathcal{H}^2_\beta).\]
\item there is a constant $\kappa>0$ such
that
$$
    \sum_{-\lambda_k \in Q_I}  \frac{|2\re  \lambda_k|^{2}|1-\lambda_k |^{2\beta}}{|b_{\infty,k}|^2
    |g_k|^{2} } \le \kappa |I|     \text{ for all
    intervals } I \subset i\RR.
$$
\end{enumerate}
\end{theorem}

Similar results can be proved for null-controllability, along the lines of \cite{jp06}.

\section*{Acknowledgements}
This work was supported by EPSRC grant EP/I01621X/1.
The third author also acknowledges partial support of this work by a Heisenberg Fellowship of the German Research Foundation (DFG).


\begin{thebibliography}{99}


\bibitem{DP94}
N. Das and J.R. Partington,  
Little Hankel operators on the half-plane. 
{\em Integral Equations Operator Theory\/} 20 (1994), no. 3, 306--324. 


\bibitem{DGM}
P. Duren, E.A. Gallardo-Guti\'errez and A. Montes-Rodr\'\i guez, 
A Paley--Wiener theorem for Bergman spaces with application to invariant subspaces. 
{\em Bull. Lond. Math. Soc.} 39 (2007), no. 3, 459--466. 


\bibitem{haak} B. H. Haak,  On the Carleson measure criterion in linear systems theory. 
{\em Complex Anal. Oper. Theory\/} 4 (2010), no. 2, 281--299.
  
\bibitem{zen09}
Z. Harper,   Boundedness of convolution operators and input-output maps between weighted spaces. 
{\em Complex Anal. Oper. Theory\/} 3 (2009), no. 1, 113--146. 

\bibitem{zen10}
Z. Harper,  Laplace transform representations and Paley--Wiener theorems for functions on vertical strips. 
{\em Doc. Math.} 15 (2010), 235--254.
  

\bibitem{HR83}
L.F.~Ho and D.L.~Russell,
Admissible input elements for systems in Hilbert space and a Carleson measure criterion,
{\em SIAM J. Control Optim.}  21  (1983),  no. 4, 614--640.  Erratum,
{\em SIAM J. Control Optim.}  21  (1983),  no. 6, 985--986.

\bibitem{jp}
B.~Jacob and J.R.~Partington, 
Admissibility of control and observation operators for semigroups: a survey.  
{\em Current trends in operator theory and its applications}, 199--221, 
Oper. Theory Adv. Appl., 149, Birkh\"auser, Basel, 2004. 


\bibitem{jp06}
B.~Jacob and J.R.~Partington, 
On controllability of diagonal systems with one-dimensional input space.  
{\em Systems Control Lett.}  55  (2006),  no. 4, 321--328. 

\bibitem{jpp07}
B. Jacob, J.R. Partington and S. Pott, Interpolation by vector-valued analytic functions, with applications to controllability. 
{\em  J. Funct. Anal.}  252  (2007),  no. 2, 517--549. 

\bibitem{JPP09}
B. Jacob, J.R. Partington and S. Pott, Tangential interpolation in weighted vector-valued $H^p$ spaces, with applications. 
{\em Complex Anal. Oper. Theory\/} 3 (2009), no. 3, 697--727. 

\bibitem{jpp12}
B. Jacob, J.R. Partington and S. Pott,
On Laplace--Carleson embedding theorems, submitted 2011.
{\tt http://arxiv.org/abs/1201.1021}




\bibitem{mcphail}
J.D.~McPhail,
A weighted interpolation problem for analytic functions.
{\em Studia Math.} 96 (1990), no. 2, 105--116.




\bibitem{ss}
H.S.~Shapiro and A.L.~Shields,
On some interpolation problems for analytic functions.
{\em Amer. J. Math.} 83 (1961), 513--532.

\bibitem{staffans}
O. Staffans, {\em Well-posed linear systems}. Encyclopedia of Mathematics and its Applications, 103. Cambridge University Press, Cambridge, 2005. 

\bibitem{steg} D.A. Stegenga, Multipliers of the Dirichlet space. {\em Illinois J. Math.} 24 (1980), no. 1, 113�139
\bibitem{stein}
E.M. Stein, {\em Topics in harmonic analysis related to the Littlewood--Paley theory}. Annals of Mathematics Studies, No. 63 Princeton University Press, Princeton, N.J.; University of Tokyo Press, Tokyo, 1970.




\bibitem{TW}
M. Tucsnak and G. Weiss, {\em Observation and control for operator semigroups}. Birkh\"auser Advanced Texts: Basler Lehrb\"ucher.  Birkh\"auser Verlag, Basel, 2009. 

\bibitem{weiss88}
G. Weiss,  Admissibility of input elements for diagonal semigroups on $l\sp 2$,
{\em Systems Control Lett.} 10  (1988),  no. 1, 79--82.

\bibitem{weiss99}
G. Weiss, A powerful generalization of the Carleson measure theorem? {\em Open problems in mathematical systems and control theory}, 267--272, Comm. Control Engrg. Ser., Springer, London, 1999.

\bibitem{wynn}
A. Wynn,   $\alpha$-admissibility of observation operators in discrete and continuous time. 
{\em Complex Anal. Oper. Theory\/} 4 (2010), no. 1, 109--131. 

\bibitem{wynn1}
A. Wynn,   Counterexamples to the discrete and continuous weighted Weiss conjectures. {\em SIAM J. Control Optim.} 48 (2009), no. 4, 2620--2635


\end{thebibliography}
\end{document}